\documentclass[reqno,a4paper]{amsart}
\usepackage{amssymb,graphicx,eepic,rotating,wrapfig,psfrag}
\usepackage{array,multirow,comment,hyperref}
\usepackage[all]{xy}
\xyoption{import}
\usepackage{longtable}
\usepackage{caption}
\usepackage{booktabs}
\excludecomment{maximacode}
\excludecomment{extratable}
\excludecomment{myproof}
\includecomment{mytable}

\def\Q{{\bf Q}}
\def\Z{{\bf Z}}
\def\C{{\bf C}}

\def\smfld{{$7$-manifold\ }}
\def\smflds{{$7$-manifolds\ }}
\def\SU{{\rm SU}}
\def\U{{\rm U}}

\def\diag#1{{\rm diag}(#1)}
\def\stf#1#2{{{\rm V}_{#1}(\C^{#2})}}
\def\one{{\bf 1}}
\def\Lk{{\rm Lk}}
\def\ds{{\displaystyle}}
\def\sgn{{\rm sign}}

\def\autkr{{\mathfrak G}}
\def\eks{{Eschenburg-Kruggel space}}
\def\absr{{|r|}}
\def\cpp{{\tt C++}}
\def\longint{{\tt signed long int}}
\def\longintminmax{{\tt signed long int} $[-2\,147\,483\,648,2\,147\,483\,647]$}
\def\longintminmax{the set of {\tt signed long int}s is $[-2 \times 1024^3,2 \times 1024^3-1]$}
\def\longintminmax{the set of {\tt signed long int}s equals $[-2^{31},2^{31}-1] \cap \Z$}
\def\gnu{{\sc Gnu}}
\def\gnugmp{{\sc Gnu Gmp}}
\def\gmp{{\sc Gmp}}
\def\gmpfrxx{{\sc Gmpfrxx}}
\def\maple{{\sc Maple}}
\def\maxima{{\sc Maxima}}
\def\bc{{\sc Bc}}
\def\urltilda{\kern -.15em\lower .7ex\hbox{\~{ }}\kern .04em}
\def\urldot{\kern -.10em.\kern -.10em}
\def\urlhttp{http\kern -.10em\lower -.1ex\hbox{:}\kern -.12em\lower 0ex\hbox{/}\kern -.18em\lower 0ex\hbox{/}}
\def\mywebpageurl#1{{http://www.maths.ed.ac.uk/~lbutler/#1}}
\def\mywebpagetxt#1{{http://www.maths.ed.ac.uk/\urltilda{}lbutler/#1}}
\def\minhcheadera{{Smooth structures on a}}
\def\minhcheaderb{{homeomorphism class of \eks{s}}}
\def\minhcheaderc#1#2#3#4{{with $|r|=#1$, $s=#2$, $p_1=#3$, $s_2=#4$}}
\newcounter{thmctr}
\newcounter{defnctr}
\newcounter{rmkctr}
\newcounter{propctr}
\newcounter{quctr}
\newtheorem{thm}[thmctr]{Theorem}
\newtheorem{defn}[defnctr]{Definition}
\newtheorem{remark}[rmkctr]{Remark}
\newtheorem{proposition}[propctr]{Proposition}
\newtheorem{question}[quctr]{Question}

\begin{document}

\title{Smooth structures on Eschenburg spaces: numerical computations}
\address{School of Mathematics and the Maxwell Institute of
  Mathematical Sciences, 6214 James Clerk Maxwell Building,
  The University of Edinburgh,  Edinburgh, UK, EH9 3JZ}
\email{l.butler@ed.ac.uk}
\urladdr{http://www.maths.ed.ac.uk/~lbutler/}
\author[Butler]{Leo T. Butler}
\date{\today}

\thanks{The author thanks T. Chinburg, C. Escher and W. Ziller for
  helpful comments and providing their \maple{} code; J. Wilkening for
  releasing his \gmpfrxx{} code; M. Kreck for discussions of
  \cite{Kr1}; and C. Peterson.}

\begin{abstract}
This paper numerically computes the topological and smooth invariants
of \eks{s} with small fourth cohomology group, following Kruggel's
determination of the Kreck-Stolz invariants of special Eschenburg
spaces~\cite{Esch,Kr1,KS1}. The \gnugmp{} arbitrary-precision
library is utilised~\cite{gmp}.
\end{abstract}


\maketitle

\section{Introduction} \label{sec:intro}

In \cite{AW}, Aloff and Wallach introduced a family of \smflds that
are homogeneous spaces of $\SU_3$ as follows: let $p,q$ be coprime
integers and let $U_{p,q} \subset \SU_3$ be the subgroup of diagonal
matrices of the form $\diag{z^p,z^q,z^{-p-q}}$ for $z\in S^1$. The
Aloff-Wallach \smfld $M_{pq} = \SU_3/U_{p,q}$. Aloff and Wallach
showed that a bi-invariant metric on $\SU_3$ induced a
positively-curved submersion metric on the quotient $M_{p,q}$. In
\cite{KS1,KS2}, Kreck and Stolz studied the topological and smooth
classification of Aloff-Wallach spaces. Amongst other things, they
showed that there are diffeomorphic Aloff-Wallach spaces that are not
$\SU_3$ equivariantly diffeomorphic: the `smallest' example occurs
with $(p,q)$ equal to $(-4\,638\,661,582\,656)$ and
$(-2\,594\,149,5\,052\,965)$! \cite[p. 468]{KS2} Each of these spaces
has a finite cyclic fourth integral cohomology group; they showed,
through a computer search, that if the order of $H^4(M_{p,q};\Z)$ is
less than $r=2\,955\,27\,597$, then the topological structure
determines the smooth structure. Additional computer searches,
attributed to Zagier and Odlyzko, revealed homeomorphic, but not
diffeomorphic, Aloff-Wallach spaces with the rank of $H^4(M_{p,q};\Z)$
between the above number and roughly $2\times 10^{20}$
\cite[p. 467]{KS2}. In all cases, there were no reported examples of a
topological Aloff-Wallach space whose $28$ distinct smooth structures
are themselves diffeomorphic to Aloff-Wallach spaces.

In \cite{Esch}, Eschenburg introduced a family of \smflds that
generalise Aloff-Wallach spaces \cite{Kr1}. Let $U \cong \U_1$ be a
subgroup of $\U_3 \times \U_3$ such that the natural action of $U$ on
$\U_3$ defined by
\begin{align} \label{al:u}
\forall u=(u_1,u_2) \in U, g \in \U_3:&& u\cdot g = u_1 g u_2^{-1}
\end{align}
stabilises $\SU_3$ and is free. $U$ is conjugate to a diagonal
subgroup characterised by $2$ integer vectors $k$ and $l\in\Z^3$ such
that $k_0+k_1+k_2=l_0+l_1+l_2$:
\begin{align} \label{al:kl}
U_{kl} &= \left\{  \diag{z^{k_0},z^{k_1},z^{k_2}} \oplus \diag{z^{l_0},z^{l_1},z^{l_2}}\
:\ z \in S^1\right\}.
\end{align}
The freeness of the action \eqref{al:u} is equivalent to the property
that
\begin{align} \label{al:free}
\forall \text{ permutations } \sigma: && k-\sigma(l) \text{ is a
  primitive vector in } \Z^3.
\end{align}
Eschenburg defined $k,l$ to be {\em admissible} if
\begin{align} \label{al:gcd}
\begin{array}{lllc}
\gcd(k_0-l_0,k_1-l_1), & \gcd(k_0-l_0,k_1-l_2), \\
\gcd(k_0-l_1,k_1-l_0), & \gcd(k_0-l_1,k_1-l_2),  \\
\gcd(k_0-l_2,k_1-l_0), & \textrm{and } \gcd(k_0-l_2,k_1-l_1) & \textrm{equal } 1.\\
\end{array}
\end{align}

\begin{defn} \label{de:espace}
Let $k,l \in \Z^3$ satisfy $k_0+k_1+k_2=l_0+l_1+l_2$ and the
admissibility conditions \eqref{al:gcd} and define $U_{kl}$ as in
\eqref{al:kl}. The $7$-manifold $E_{kl} := \SU_3 / U_{kl}$ is called
an {\em Eschenburg} space.
\end{defn}

Eschenburg computed the integral cohomology ring of
$E_{k,l}=\SU_3/U_{k,l}$ and proved that these spaces were strongly
inhomogeneous in most cases. He also showed that under certain
conditions on $k,l$, a bi-invariant metric on $\SU_3$ induces a
positively-curved submersion metric on $E_{k,l}$.

In \cite{Kr1}, Kruggel computed the Kreck-Stolz invariants of a broad
number of Eschenburg spaces---henceforth an \eks---and obtained a
classification of these \eks{s} up to homotopy, homeomorphism and
diffeomorphism. In \cite{CEZ}, Chinburg, Escher and Ziller implemented
a computer search for homeomorphic, but not diffeomorphic,
positively-curved (resp. 3-Sasakian) \eks{s}. They found that for $\#
H^4(E_{k,l};\Z) < 8000$, there is a unique pair of homeomorphic, but
not diffeomorphic, positively-curved \eks{s}. In \cite{Bu}, the
present author proved that the existence of a real-analytically
completely integrable convex Hamiltonian is a non-trivial smooth
invariant of the configuration space, and proved the complete
integrability of geodesic flows on all \eks{s}. That work motivated
the

\begin{question} \label{qu:smooth-structures}
Let $E$ be a topological Eschenburg space. Is each smooth structure on
$E$ diffeomorphic to an Eschenburg space $E_{k,l}$?
\end{question}

One knows from the work of Kreck and Stolz that each topological
Eschenburg space admits $28$ distinct oriented smooth structures, but
one does not know if each structure is represented by an Eschenburg
space. From the above-mentioned results, it is not clear whether each
distinct oriented smooth structure on a topological Eschenburg space
is represented by an Eschenburg space or if such representatives are
rather sparse, as for Aloff-Wallach spaces. This note attempts to cast
some light on this question.

\begin{thm} \label{th:1}
Let $I=[-850,850]$ and $J=[1,101]$. Amongst the \eks{s} with $(k,l)
\in I^3 \times I^3$, for each odd $\absr = \# H^4(E;\Z)$ in the
interval $J$, columns 2 \& 9 of Table \ref{ta:hc} show a lower bound
on the number of oriented homeomorphism classes. For $|r| \leq 9$,
each oriented homeomorphism class of \eks{s} has each of its $28$
distinct oriented smooth structures represented by an \eks{} $E_{k,l}$
with $(k,l) \in I^3 \times I^3$.
\end{thm}

\begin{remark} \label{re:1}
{\rm 

Columns 3--7 \& 10--14 of Table \ref{ta:hc} list the number of
topological \eks{s}, for a fixed $|r|$, that have the stated number of
oriented smooth structures represented by \eks{s}.

The smooth structures on a topological \eks{} are an orbit of the
group of homotopy $7$-spheres ($\cong \Z_{28}$). The Kreck-Stolz
invariant $s_1$ is additive under this action: if $\Sigma$ is a
homotopy $7$-sphere and $E$ is an \eks{}, then
$s_1(E\#\Sigma)=s_1(E)+s_1(\Sigma)$ and $28 \cdot s_1(\Sigma) \equiv 0
\bmod 1$. This implies that each topological \eks{} has $28$ distinct
oriented smooth structures \cite{Kr1}. The difficulty is that the
surgery description of the smooth structure $E \# \Sigma$ does not
appear to contain information about the structure of $E \# \Sigma$ as
a \eks{}.

Tables \ref{ta:eks-r-1}--\ref{ta:eks-r-5} list representative \eks{s}
for each smooth structure on each topological \eks{} with $\absr \leq
5$ that was found in constructing table \ref{ta:hc}. It seems likely
that all topological and smooth \eks{s} with $\absr \leq 5$ are
enumerated in tables \ref{ta:eks-r-1}--\ref{ta:eks-r-5}.

}
\end{remark}

\begin{mytable}
{
\def\small{\tiny}
\tiny
{
\def\Z{{\bf Z}}
\def\rk{{\rm rank\,}}
\def\dash{${\small -}$}
\def\d{{}}
\newcolumntype{C}{>{$}r<{$}}
\begin{longtable}[c]{C*{1}{C}*{5}{C}|C*{1}{C}*{5}{C}}
\caption{\small $|r| = \rk H^4(E;\Z)$ versus the number of homeomorphism classes (\#Top.), and the number of homeomorphism classes with the $n$ smooth structures represented by \eks{s}, for $n=28,27$, $14 \leq n \leq 26$, $2  \leq n \leq 13$ and $n=1$.}  \label{ta:hc}
\\\toprule
\multicolumn{2}{c}{} & \multicolumn{5}{c|}{Counts} & \multicolumn{2}{c}{} & \multicolumn{5}{c}{Counts} \\
|r| & \#Top. & 28 & 27 & 14\dash26 & 2\dash13 & 1 & |r| & \#Top. & 28 & 27 & 14\dash26 & 2\dash13 & 1 \\\midrule
\endfirsthead
\toprule
\multicolumn{14}{c}{\small\slshape Table \ref{ta:hc}, continued from previous page} \\
\multicolumn{2}{c}{} & \multicolumn{5}{c}{Counts} & \multicolumn{2}{c}{} & \multicolumn{5}{c}{Counts} \\
|r| & \#Top. & 28 & 27 & 14\dash26 & 2\dash13 & 1 & |r| & \#Top. & 28 & 27 & 14\dash26 & 2\dash13 & 1 \\\midrule
\endhead
\multicolumn{14}{|c|}{\small\slshape continued next page} \\\bottomrule
\endfoot
\bottomrule
\endlastfoot
1	&	12	&	12 & \d & \d & \d & \d 	&	3	&	8	&	8 & \d & \d & \d & \d \\
5	&	48	&	48 & \d & \d & \d & \d 	&	7	&	120	&	120 & \d & \d & \d & \d \\
9	&	24	&	24 & \d & \d & \d & \d 	&	11	&	360	&	354 & 4 & 2 & \d & \d \\
13	&	576	&	542 & 22 & 12 & \d & \d 	&	15	&	32	&	32 & \d & \d & \d & \d \\
17	&	1152	&	988 & 68 & 96 & \d & \d 	&	19	&	1512	&	1216 & 86 & 204 & 6 & \d \\
21	&	80	&	64 & 10 & 6 & \d & \d 	&	23	&	2640	&	1726 & 276 & 598 & 40 & \d \\
25	&	240	&	240 & \d & \d & \d & \d 	&	27	&	72	&	72 & \d & \d & \d & \d \\
29	&	4704	&	1656 & 814 & 2212 & 22 & \d 	&	31	&	5760	&	1506 & 794 & 3080 & 380 & \d \\
33	&	240	&	114 & 30 & 92 & 4 & \d 	&	35	&	480	&	230 & 90 & 160 & \d & \d \\
37	&	8634	&	904 & 918 & 5728 & 1072 & 12 	&	39	&	384	&	118 & 58 & 176 & 32 & \d \\
41	&	11988	&	376 & 636 & 8778 & 2176 & 22 	&	43	&	12600	&	272 & 500 & 9412 & 2414 & 2 \\
45	&	96	&	60 & 20 & 16 & \d & \d 	&	47	&	17108	&	82 & 248 & 10950 & 5812 & 16 \\
49	&	1848	&	1028 & 310 & 510 & \d & \d 	&	51	&	768	&	44 & 26 & 522 & 176 & \d \\
53	&	22456	&	46 & 122 & 11414 & 10836 & 38 	&	55	&	1440	&	320 & 200 & 752 & 168 & \d \\
57	&	1008	&	28 & 36 & 666 & 278 & \d 	&	59	&	29902	&	10 & 32 & 10662 & 19034 & 164 \\
61	&	32874	&	22 & 76 & 9468 & 22764 & 544 	&	63	&	240	&	60 & 10 & 164 & 6 & \d \\
65	&	2304	&	178 & 220 & 1332 & 574 & \d 	&	67	&	39854	&	12 & 28 & 7890 & 31108 & 816 \\
69	&	1756	&	\d & \d & 864 & 878 & 14 	&	71	&	47544	&	\d & \d & 6596 & 39738 & 1210 \\
73	&	48034	&	2 & 10 & 6090 & 40864 & 1068 	&	75	&	160	&	50 & 42 & 68 & \d & \d \\
77	&	3600	&	332 & 112 & 1914 & 1234 & 8 	&	79	&	59046	&	\d & 4 & 4508 & 51962 & 2572 \\
81	&	216	&	188 & 22 & 6 & \d & \d 	&	83	&	67340	&	\d & \d & 3544 & 59816 & 3980 \\
85	&	4602	&	28 & 82 & 2670 & 1800 & 22 	&	87	&	3128	&	\d & \d & 580 & 2522 & 26 \\
89	&	78944	&	\d & \d & 2068 & 70090 & 6786 	&	91	&	5740	&	256 & 154 & 2502 & 2734 & 94 \\
93	&	3788	&	\d & \d & 468 & 3182 & 138 	&	95	&	6016	&	18 & 22 & 2836 & 3054 & 86 \\
97	&	91772	&	\d & \d & 1484 & 79690 & 10598 	&	99	&	720	&	12 & 40 & 474 & 194 & \d \\
101	&	100490	&	\d & \d & 742 & 87290 & 12458 	&	&&&&&&
\end{longtable}
}

\begin{maximacode}
tots : [[1,12,12,0,0,0,0],[3,8,8,0,0,0,0],
[5,48,48,0,0,0,0],[7,120,120,0,0,0,0],
[9,24,24,0,0,0,0],[11,360,354,4,2,0,0],
[13,576,542,22,12,0,0],[15,32,32,0,0,0,0],
[17,1152,988,68,96,0,0],[19,1512,1216,86,204,6,0],
[21,80,64,10,6,0,0],[23,2640,1726,276,598,40,0],
[25,240,240,0,0,0,0],[27,72,72,0,0,0,0],
[29,4704,1656,814,2212,22,0],[31,5760,1506,794,3080,380,0],
[33,240,114,30,92,4,0],[35,480,230,90,160,0,0],
[37,8634,904,918,5728,1072,12],[39,384,118,58,176,32,0],
[41,11988,376,636,8778,2176,22],[43,12600,272,500,9412,2414,2],
[45,96,60,20,16,0,0],[47,17108,82,248,10950,5812,16],
[49,1848,1028,310,510,0,0],[51,768,44,26,522,176,0],
[53,22456,46,122,11414,10836,38],[55,1440,320,200,752,168,0],
[57,1008,28,36,666,278,0],[59,29902,10,32,10662,19034,164],
[61,32874,22,76,9468,22764,544],[63,240,60,10,164,6,0],
[65,2304,178,220,1332,574,0],[67,39854,12,28,7890,31108,816],
[69,1756,0,0,864,878,14],[71,47544,0,0,6596,39738,1210],
[73,48034,2,10,6090,40864,1068],[75,160,50,42,68,0,0],
[77,3600,332,112,1914,1234,8],[79,59046,0,4,4508,51962,2572],
[81,216,188,22,6,0,0],[83,67340,0,0,3544,59816,3980],
[85,4602,28,82,2670,1800,22],[87,3128,0,0,580,2522,26],
[89,78944,0,0,2068,70090,6786],[91,5740,256,154,2502,2734,94],
[93,3788,0,0,468,3182,138],[95,6016,18,22,2836,3054,86],
[97,91772,0,0,1484,79690,10598],[99,720,12,40,474,194,0],
[101,100490,0,0,742,87290,12458]];
map(lambda([t],t[2]-sum(t[i],i,3,7)),tots);
\end{maximacode}

}
\end{mytable}

This note is structured as follows: section \ref{sec:eks} reviews
Kruggel's condition C; section \ref{sec:kr} reviews Kruggel's
computation of the Kreck-Stolz invariants; section \ref{sec:method}
explains how the Kreck-Stolz invariants were computed in software; and
appendices \ref{ssec:app-a}--\ref{ssec:app-d} add several tables.

\section{\eks{s}} \label{sec:eks}

To compute the Kreck-Stolz invariants of Eschenburg spaces, Kruggel
observed that the projection of an $x \in \SU_3$ onto its first two
columns in the Steifel manifold $\stf{2}{3}$ is a diffeomorphism. From
the embedding of $\stf{2}{3} \subset\C^{2\cdot3}$, Kruggel constructed
an $8$-manifold $W'$ with boundary $\stf{2}{3}$. The action of
$U_{kl}$ descends naturally to $\C^{2\cdot 3}$ and $W'$, but the
action on $W'$ has $3$ singular orbits. One can cut away these three
singular orbits to construct a cobordism between $E_{kl}$ and a union
of $3$ lens spaces -- provided that the matrix
\begin{align}  \label{al:C}
A &=
\begin{bmatrix}
k_0-l_0 & k_0-l_1 & k_0-l_2\\
k_1-l_0 & k_1-l_1 & k_1-l_2\\
k_2-l_0 & k_2-l_1 & k_2-l_2
\end{bmatrix}
\end{align}
has a column or row containing non-zero, pairwise coprime entries.

\begin{defn}[Kruggel 2006] \label{de:eks}
The Eschenburg space $E_{k,l}$ satisfies {\em condition C} iff the
matrix $A$ has a column or row containing non-zero, pairwise coprime
entries.  An Eschenburg space that satisfies condition C is called an
{\em Eschenburg-Kruggel space}.
\end{defn}

\begin{remark} \label{re:2}
{\rm

Note that the coprimality conditions \eqref{al:gcd} do not imply that
all entries of $A$ are non-zero. The Eschenburg space $E_{kl}$ with
$k=(-1,-1,2)$ and $l=(-2,0,2)$ has
\begin{align} \label{al:ce}
A&=
\begin{bmatrix}
1 & -1 & -3\\
1 & -1 & -3\\
4 & 2  & 0
\end{bmatrix}.
\end{align}
This defines an Eschenburg-Kruggel space according to definition
\ref{de:eks}. Indeed, the coprimality conditions \eqref{al:gcd} are
satisfied, since they are
\begin{align} \label{al:gcd-a}
\begin{array}{lllcll}
\gcd(A_{00},A_{11}),& \gcd(A_{00},A_{12}),&&      \gcd(1,-1), & \gcd(1,-3),\\
\gcd(A_{01},A_{10}),& \gcd(A_{01},A_{12}),&i.e.&  \gcd(-1,1), & \gcd(-1,-3),\\
\gcd(A_{02},A_{10}),& \gcd(A_{02},A_{11}) &&      \gcd(-3,1), & \textrm{and } \gcd(-3,-1)
\end{array}
\end{align}
which are all unity; and condition C is satisfied by the left-most
column of $A$. See remark \ref{re:3} for more.

}
\end{remark}

\section{Invariants of Eschenburg-Kruggel spaces} \label{sec:kr}
Let $E_{kl}$ be an Eschenburg space. Let $u$ be the Chern class of the
bundle $S^1=U_{kl} \hookrightarrow \SU_3 \to E_{kl}$. Eschenburg
proved that the non-trivial parts of the integral cohomology ring of
$E_{k,l}$ has the following structure:
\begin{align} \label{al:coho}
H^2(E_{kl};\Z)&= \Z\cdot u, & H^4(E_{kl};\Z)&= \Z_{r}\cdot u^2.
\end{align}
The integer $r=\sigma_2(k)-\sigma_2(l)$ where $\sigma_j$ is the
$j$-th elementary symmetric polynomial, $\sigma_j(x) = \sum_{i_1 <
  \cdots < i_j} x_{i_1} \cdots x_{i_j}$. The linking form of $E_{k,l}$
is plainly determined by the linking number of $u^2$ with
itself. Kruggel showed that this equals
\begin{align} \label{al:lk}
\Lk(u^2,u^2) &= -\frac{s^{-1}}{r} \bmod 1,
\end{align}
where $s=\sigma_3(k)-\sigma_3(l)$ and $s^{-1}$ is the multiplicative
inverse of $s \bmod r$. Kruggel also showed that the first
Pontryagin class of $E_{kl}$ equals
\begin{align} \label{al:p1}
p_1(E_{kl}) &= p_1\cdot u^2 \bmod r &&\text{where } p_1= 2 \sigma_1(k)^2 - 6
\sigma_2(k).
\end{align}
Although this expression appears to be asymmetric in $k$ and $l$, the
sum condition plus the definition of $r$ ensures that it is
well-defined.

In addition to the above invariants, Kruggel was able to compute the
Kreck-Stolz invariants for Eschenburg-Kruggel spaces. To explain, let
$p\neq 0$ be coprime to the non-zero integers $p_0,\ldots,p_3$, and let
\begin{align} \label{al:lens}
&L=L(p;p_0,p_1,p_2,p_3) = S^7/C \text{ where }\\
&C=\{ \diag{e^{\frac{2\pi ik p_0}{p}},e^{\frac{2\pi ik p_1}{p}},e^{\frac{2\pi ik p_2}{p}},e^{\frac{2\pi ik p_3}{p}}}\ :\ k=0,\ldots,p-1\}
\notag
\end{align}
be a lens space. Define the following functions
\begin{align} \label{al:siL}
s_1(L) &= \frac{1}{2^7\cdot7\cdot p}\ds\sum_{k=1}^{|p|-1}
\prod_{j=0}^3 \cot\left(\frac{k\pi p_j}{p}\right) + \frac{1}{2^4 \cdot
  p}\ds\sum_{k=1}^{|p|-1} \prod_{j=0}^3 \csc\left(\frac{k\pi
  p_j}{p}\right) &&\bmod 1\\
s_2(L) &= \frac{1}{2^4 \cdot
  p}\ds\sum_{k=1}^{|p|-1} \left(e^{\frac{2\pi
    ik}{p}}-1\right)\prod_{j=0}^3 \csc\left(\frac{k\pi p_j}{p}\right)
&& \bmod 1 \notag
\end{align}
These are the Kreck-Stolz invariants of the lens space $L$ in
\eqref{al:lens}, and they take values in $\Q/\Z$.

Assume that the left-most column of the matrix $A$ has pairwise
coprime, non-zero entries (the remaining cases are described
below). The above-described cobordism exhibits $E_{kl}$ as cobordant
to the disjoint union of the three lens spaces:
\begin{align} \label{al:le}
L_0=L(A_{00};A_{10},A_{20},A_{11},A_{21}) &&
L_1=L(A_{10};A_{00},A_{20},A_{01},A_{21}) \notag\\
L_2=L(A_{20};A_{00},A_{10},A_{01},A_{11}).
\end{align}
Let us see that $L_0$ is indeed a lens space. By condition C, the
integers $A_{j0}$ are pairwise coprime and non-zero. The primitivity
condition \eqref{al:free} implies that $A_{00}$ is coprime to $A_{11}$
and $A_{21}$. For example, suppose that $A_{00}$ and $A_{11}$ have a
divisor $d>1$, so one can write $A_{00}+A_{11}=dc$. The condition that
$\sum k_i=\sum l_i$ in definition \ref{de:espace} is equivalent to
$A_{00}+A_{11}+A_{22}=0$, so $A_{22}=-dc$. If $c=0$, then the vector
$k-l=(A_{00},A_{11},A_{22})$ is not primitive; if $c \neq 0$, then the
same vector is not primitive, too. (This argument also shows that if
$A_{22}=0$, then $A_{00}=-A_{11}=\pm 1$.) The remaining verifications
for $L_1$ and $L_2$ are similar.

Kruggel showed that the Kreck-Stolz invariants are equal to
\begin{align}
s_1(E_{kl}) &= \frac{\sgn(w)}{2^5\cdot7} -
\frac{q^2}{2^7\cdot7\cdot w}
 -\sum_{i=1}^3 s_1(L_i) &&\bmod 1\label{al:s1} \\
s_2(E_{kl}) &= \frac{q-2}{2^4\cdot3\cdot w}
 -\sum_{i=1}^3 s_2(L_i) &&\bmod 1\label{al:s2}
\intertext{where}
q&=A_{00}^2 + A_{10}^2 + A_{20}^2 + A_{01}^2 + A_{11}^2 + A_{21}^2 -
(l_0-l_1)^2, \label{al:q} \\
w&=r\cdot A_{00}A_{10}A_{20}. \label{al:w}
\end{align}
These invariants are transcendental functions of the variables
$k,l$. This fact, plus the fact that the sums can have a rather large
number of terms, means that showing two \eks{s} are homeomorphic or
diffeomorphic is rather difficult. However, since $s_i$ is a rational
integer, one can use a few numerical tricks to prove equality of these
invariants.

\begin{remark}[{\it c.f.} remark \ref{re:2}] \label{re:3}
{\rm

The well-definedness of Kruggel's formulae (\ref{al:s1}--\ref{al:s2})
amounts to the statement that if $E_{k,l}$ satisfies condition C, then
$w$ \eqref{al:w} does not vanish. Indeed, from the remark above, the
lens-space invariants $s_j(L_i)$ (\ref{al:siL}--\ref{al:le}) are
well-defined if $w \neq 0$. Since condition C is assumed to hold for
the left-most column of $A$, $A_{00}A_{10}A_{20} \neq 0$. In addition,
Kruggel proved that $r$ must be odd~\cite[p. 572]{Kr1} (in fact, since
$H^3(E_{k,l};\Z)$ vanishes, Poincar\'e duality implies $r \neq
0$). Therefore, $w \neq 0$.

}
\end{remark}

The results of this note rely on

\begin{thm}[Kruggel 2005] \label{th:kruggel}
Two \eks{s}, $E_{k,l}$ and $E_{k',l'}$ are orientation-preserving {\em
  homeomorphic} if $|r|,s,p_1$ and $s_2$ coincide. If, in addition,
$s_1$ coincides, then they are orientation-preserving {\em
  diffeomorphic}.
\end{thm}

\subsection{Automorphisms and invariants} \label{ssec:aut}
To compute the Kreck-Stolz invariants of \eks{s} in general, one uses
the extension of the natural action of Weyl group of $\SU_3 \times
\SU_3$ by the automorphism that interchanges factors. Concretely, let
$S_3$ be the symmetric group acting naturally on $\Z^3$ by
permutations, let $\tau$ be the involutive automorphism of $\Z^3
\oplus \Z^3$ which acts by $(k,l) \mapsto (l,k)$ and let $\eta : (k,l)
\mapsto (-k,-l)$.

The group generated by $S_3 \times S_3$, $\tau$ (resp. $S_3 \times
S_3$, $\tau$ and $\eta$) is denoted by $\autkr^+$
(resp. $\autkr$). $\autkr$ is a group of order $144$ and $\autkr^+$ is
an index $2$ subgroup.

\begin{proposition}[\cite{Esch}] \label{pr:aut}
For each $\sigma \in \autkr$, the Eschenburg spaces $E_{k,l}$ and
$E_{\sigma(k,l)}$ are diffeomorphic. If $\sigma\in\autkr^+$, they are
orientation-preserving diffeomorphic.
\end{proposition}

\begin{myproof}
Let $\sigma \in \autkr$ and let $k,l \in \Z^3$ define an Eschenburg
space $E_{k,l}$. One seeks to construct diffeomorphisms $s$
(resp. $S$) of Eschenburg spaces (resp. $\SU_3$) such that the
following diagram commutes
\begin{align} \label{al:xy}
\xymatrix@M=6pt@C+3pt{
\SU_3 \ar[d]_{/U_{k,l}} \ar[r]^{S} & \SU_3
\ar[d]^{/U_{\sigma(k,l)}}\\
E_{k,l} \ar[r]^{s} & E_{\sigma(k,l)}.
}
\end{align}
We split the verification into 3 cases:
\begin{enumerate}
\item $\sigma \in S_3 \times S_3$: Let $\Delta= \left\{ (g,h) \in \U_3
\times \U_3\ :\ \det g = \det h \right\}$. This is a connected Lie
group that acts naturally on $\SU_3$ via
\begin{align} \label{al:a}
(g,h)\cdot u &= guh^{-1} & \forall (g,h)\in \Delta, u \in \SU_3.
\end{align}
If $W$ is the Weyl group of $\SU_3$, then $W \times W \subset
\Delta$. The action by conjugation of $W$ on the maximal torus $T$ of
diagonal matrices in $\SU_3$ induces the permutation representation of
$S_3$ on $\Z^3$. In this case, the transformation $S=(s_1,s_2) \in W
\times W$ acting via \eqref{al:a} induces the transformation $s$
defined by
\begin{align}
s(g\cdot U_{k,l}) &= \left\{
s_1u^ks_1^{-1}(s_1gs_2^{-1})s_2u^ls_2^{-1} \right\} \label{al:s}\\
&= S(g)\cdot U_{\sigma(k,l)}, \notag\\
\intertext{where}
u^m &= \diag{z^{m_0},z^{m_1},z^{m_2}} &&\forall z \in S^1, m \in \Z^3.\notag
\end{align}
Since $W$ normalises $T$, the right-hand side of \eqref{al:s} is
well-defined independent of the representative $g$. Since $S$ is
homotopic to $id$ and the orientation on the respective fibres are
preserved, $s$ is an orientation-preserving diffeomorphism.

\item $\sigma=\eta$. Let $\overline{U}_{k,l}$ be $U_{k,l}$ with the
  opposite orientation; $\overline{U}_{k,l}$ is clearly equal to
  $U_{-k,-l}$. It is clear that $S=id$ induces the natural
  orientation-reversing diffeomorphism $s$ in this case.

\item $\sigma=\tau$. Let $S$ be the orientation-preserving
  diffeomorphism $g \mapsto g^{-1}$. $S$ induces an isomorphism
  $U_{k,l} \to U_{-l,-k}$ which preserves orientation; $S$ reverses
  orientation along the orbits of these groups, however. Therefore,
  $S$ induces the orientation-reversing diffeomorphism $s : E_{k,l}
  \to E_{-l,-k}$. From the previous step, $E_{k,l} \to E_{-l,-k} \to
  E_{l,k}$ is orientation preserving.
\end{enumerate}
\end{myproof}

\begin{remark} \label{re:perm}
{\rm

With the above proposition, the formulae for the Kreck-Stolz
invariants can be extended to all \eks{s} as follows. The
Eschenburg space $E_{kl}$ is orientation-preserving diffeomorphic to
$E_{\alpha(k),\beta(l)}$ for any permutations $\alpha,\beta \in S_3$. In
addition, $E_{k,l}$ is orientation-preserving diffeomorphic to
$E_{l,k}$. The permutation $\alpha$ permutes the rows (resp. $\beta$
permutes the columns) of $A$, while the diffeomorphism $E_{kl}\to
E_{lk}$ induces $A \mapsto -A'$.

It follows that if the column $j$ (resp. row $j$) of $A$ has non-zero
pairwise coprime entries, then the leftmost column of $A_{k,\beta(l)}$
(resp. $A_{l,\beta(k)}$) has non-zero pairwise coprime entries and
$E_{kl}$ is orientation-preserving diffeomorphic to $E_{k,\beta(l)}$
(resp. $E_{l,\beta(k)}$) where $\beta=(0\, j)$. By this observation,
one can compute the Kreck-Stolz invariants of any \eks\ by means
of the formulae (\ref{al:s1},\ref{al:s2}).

The proposition also implies that each Eschenburg space $E_{k,l}$ has
a representative, up to orientation, where $k_0 \leq k_1 \leq k_2$,
$l_0 \leq l_1 \leq l_2$ and $k_0 \leq l_0$.

}
\end{remark}

\section{Methodology} \label{sec:method}
The search for homeomorphic smooth \eks{} neatly divides into three
separate searches:
\begin{enumerate}
\item search over a domain of parameters $(k,l) \in \Z^3 \times \Z^3$
  for \eks{s};
\item computation of the invariants $r$, $s$, $p_1$ and $s_1, s_2$ in
  terms of the parameters $(k,l)$;
\item search the data generated for matching invariants.
\end{enumerate}

Due to the size of the sample space considered, it was decided to do
the first two steps in compiled code. The structure of the problem
led to the choice of \cpp{} as the language of choice. 

The computations to generate all of the tables in this note took
approximately six weeks of continuous cpu time on a single core of a
2-core 3.0GHz Intel Core Duo E6850 cpu with 4MB cache and 3.3GB DDRAM
4.0GB swap. The operating system was RHEL with the 2.6.8 Linux kernel.



\subsection{The search over parameter space} \label{ssec:step-1}
Let us define the parameter space and explain how the search is
conducted.

\subsubsection{The parameter space} \label{sssec:step-1-domain}
Let $\one \in \Z^3$ be the vector whose elements are all unity. If
$E_{k,l}$ is an Eschenburg space, then $E_{k+n\one,l+n\one}$ is the
same Eschenburg space for any $n \in \Z$. There is, therefore, a
unique representative of $(k,l) + \Z(\one,\one)$ such that $\sum k_i =
\sum l_i \in [0,2]$. All searches were conducted with this
constraint.
\footnote{In tables \ref{ta:eks-r-1}--\ref{ta:eks-r-5}, one finds the
  sums reported lie in $[-2,2]$. Those spaces with sum reported in
  $[-2,-1]$ are obtained by reversing the orientation of a space whose
  sum lies in $[1,2]$.}

\subsubsection{The search} \label{sssec:step-1-search}
The speed of the arithmetic in the native \longint{} class of integers
in \cpp{} argued in favour of performing testing the admissibility
condition \eqref{al:gcd} and condition C (definition \ref{de:eks}) in
\longint. 

The coprimality tests are conducted by a two-part process. First, an
$N \times N$ lookup table is created. The $(i,j)$ entry of the lookup
table equals $1$ if $i$ and $j$ are coprime and $ij \neq 0$;
otherwise, it is $0$. If $|i|$ or $|j|$ exceed $N$, the Euclidean
algorithm is first employed to reduce both $i$ and $j$ until the
lookup table can be used. The parameter $N$ is chosen at compile time;
in our tests $N=2000$ was chosen so that all coprimality tests
required only a lookup.

\subsection{Computation of the invariants} \label{ssec:step-2} This is
broken into two parts.
\subsubsection{Integer invariants} \label{sssec:step-2-integer-invariants} 
If $(k,l)$ define an \eks{}, then the rank of $H^4(E;\Z)$, $|r|$, and
the first Pontryagin class $p_1$ were computed using \longint{}
arithmetic. Since \longintminmax{}, and both $r$ and $p_1$ are
quadratic forms in $(k,l)$, \longint{} arithmetic does not run into
under/overflow errors for $|k_i|<10922$. For the purposes of this
note, all computations of $r$ and $p_1$ were done in \longint{}
arithmetic.

Since $s$ is cubic in $(k,l)$, under/overflow does not affect
computation for $|k_i|,|l_i| < 1023$. This relatively small bound led
us to use \gmp{} arbitrary precision floats to compute $s$ (see
below).

\subsubsection{Rational invariants} \label{sssec:step-2-rational-invariants} 
From the definition of the Kreck-Stolz invariants \eqref{al:siL}, one
can see that individual terms in each summand can be $O(1/p^3)$. 

The \gnugmp{} package, along with its \gmpfrxx{} front-end for \cpp{},
permit one to do arbitrary precision arithmetic from within
\cpp{}~\cite{gmp,gmpfrxx}. Since \gnugmp{} can compute the
trigonometric functions to arbitrary precision, we elected to use this
package to compute the Kreck-Stolz invariants of an \eks{}.

The relative slowness of software-implemented arithmetic also
indicated a need to permit computation with machine-native floating
point arithmetic. The template facility of \cpp{} made it possible to
use the same code for both machine-native and software-implemented
floating-point arithmetic and allow the user to choose the precision
at run-time rather than compile-time.

\subsection{Matching invariants} \label{ssec:step-3-search} 
The final step was to match the topological and smooth invariants that
are computed for different \eks{s}. This was accomplished, in essence,
by multiple sorts. In the first step, a \cpp{} programme computed and
sorted approximately 2GB of the polynomial \eks{} invariants
($r,s,p_1$). These data were stored in text files, and these were
sorted and split according to the value of $|r|$. The Kreck-Stolz
invariants of these spaces were computed with 130 bits of precision
and stored in a second database. The resulting data were imported into
a second \cpp{} programme where homeomorphism and diffeomorphism
classes were computed. The data for tables
\ref{ta:hc}--\ref{ta:eks-r-5} were generated in this way.

\subsubsection{Testing} \label{sssec:step-4-testing}
To ensure the accuracy of the computations, several tests were
designed. These included:
\begin{enumerate}
\item replication of each of the published computations in
  \cite[section 4]{asteyetal}, table 1 of \cite{Kr1} and tables 1-6 of
  \cite{CEZ}
\footnote{In replicating these results, differing conventions for the
  projection map $x \mapsto \bar{x} \in (-\frac{1}{2},\frac{1}{2}]$
became apparent. The Chinburg-Escher-Ziller code uses the convention
that $x$ is reduced mod $1$, then $[0,\frac{1}{2}]$ is mapped to
itself by the identity and $(\frac{1}{2},1]$ is mapped to $(-1,0]$ by
    a constant shift. In our \cpp{} code, $x$ is reduced mod $1$, then
    shifted by
    $-\frac{1}{2}$.}
      \ and table 1 of \cite{Kr1};
\item replication, up to a numerical $\epsilon \sim 2^{-130}$, of
  closed-form answers for the invariants of some \eks{};
\item replication, up to a numerical $\epsilon \sim 2^{-130}$, of the
  \cpp{} computed results in \maple{}, \maxima{} and
  \bc{}~\cite{maple,maxima,bc}.
\end{enumerate}

\section{Appendices} \label{sec:app}

\subsection{Appendix A} \label{ssec:app-a}
The graph in figure \ref{gr:nvk} graphs the number $N$ of \eks{s} in
the cube $[-k,k]^6$, as a function of $k$, with the constraint that
$\sum k_i=\sum l_i \in [0,2]$. A rough heuristic indicates that
$N=O(k^4)$ for large $k$ and $\Delta N=O(k^3)$--which is nicely
captured here. It is also apparent that $\Delta N(k)$ grows like
$c_{\pm} k^3$, where $c_{\pm}$ depends only on the parity of $k$.

{
\def\Diamond{{{\cdot}}}
\def\Box{{*}}
\begin{figure}[!htb]
\input{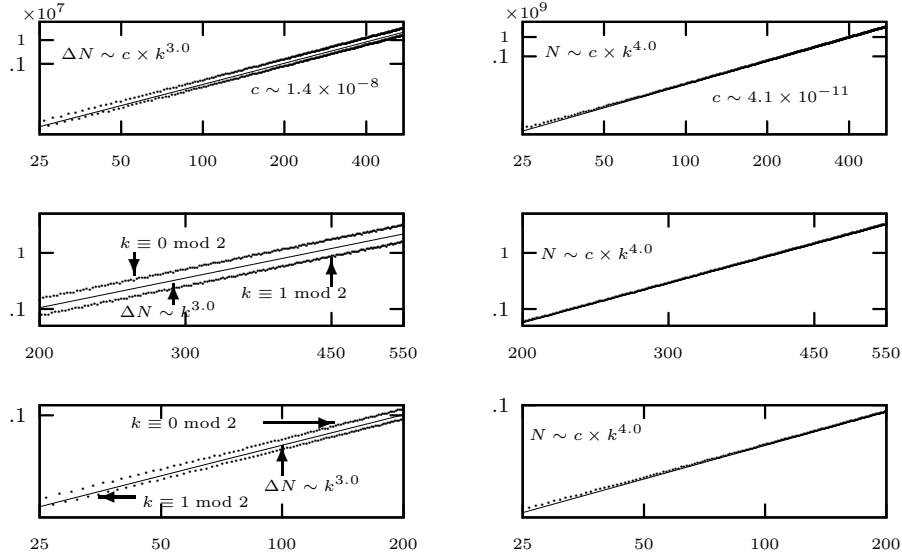}
\caption{[Log-Log scale]. The number of \eks{s}, $N=N(k)$, and the marginal number,
  $\Delta N = \Delta N(k)$, in the cube $[-k,k]^6$. Left column,
  descending: the marginal number for $k$ in the intervals $[25,550]$,
  $[200,550]$ and $[25,200]$; Right column, descending: the total
  number for the same intervals. A least-squares regression line is
  also displayed on each graph.}\label{gr:nvk}
\end{figure}
}

\subsection{Appendix B} \label{ssec:app-b}
We observed several unexplained phenomena. For fixed invariants $r,s$
and $p_1$, the Kreck-Stolz invariant $s_2$ appears to lie in the orbit
of $\Z_{n}$ acting by $x \mapsto x+\frac{1}{n} \bmod 1$ where $n=4$ or
$12$. We also observed that the values taken on by $s_1$ appear to
depend only on $|r|$, $s$ and $p_1$. 

The first columns of Tables \ref{ta:hc_1_m0_m0} and
\ref{ta:hc_3_p1_m0} shows these group actions on the Kreck-Stolz
invariants, when $|r|=1, 3$. Table \ref{ta:eks-invariants} abstracts
the picture from tables \ref{ta:hc_1_m0_m0} and \ref{ta:hc_3_p1_m0},
and shows the group actions on $s_2$ and $s_1$. It appears that
$\Z_{12}$ acts effectively except when $r \equiv 0 \bmod 3$, $r
\not\equiv 0 \bmod 3^2$, in which case $\Z_4$ acts effectively.

Additional tables are available at
\href{\mywebpageurl{kspace/}}{\mywebpagetxt{kspace/}}.

\begin{mytable}
{\tiny

{
\def\mycaption{%
  \begin{minipage}[h]{10cm}
{\tiny $\Delta s_2 = 1/12$. The values reported are convergents of
  the continued fraction approximation to the floating point
  invariants of a representative \eks{} in each class. }
\end{minipage}
}

\begin{sidewaystable}
\newcolumntype{C}{>{\tiny $}r<{$}}
\begin{tabular}{*{15}{C}} \toprule
\multicolumn{15}{c}{Computed Kreck-Stolz invariants of \eks{s}}\\

\multicolumn{1}{c}{$s_2$} & \multicolumn{14}{c}{$s_1$} \\\midrule
\multirow{2}{8mm}{$-1/2$} &-1/2&	-13/28&	-3/7&	-11/28&	-5/14&	-9/28&	-2/7&	-1/4&	-3/14&	-5/28&	-1/7&	-3/28&	-1/14&	-1/28\\%
&
0/1&	1/28&	1/14&	3/28&	1/7&	5/28&	3/14&	1/4&	2/7&	9/28&	5/14&	11/28&	3/7&	13/28\\

\multirow{2}{8mm}{$-5/12$} &-1/2&	-13/28&	-3/7&	-11/28&	-5/14&	-9/28&	-2/7&	-1/4&	-3/14&	-5/28&	-1/7&	-3/28&	-1/14&	-1/28\\%
&
0/1&	1/28&	1/14&	3/28&	1/7&	5/28&	3/14&	1/4&	2/7&	9/28&	5/14&	11/28&	3/7&	13/28\\

\multirow{2}{8mm}{$-1/3$} &-1/2&	-13/28&	-3/7&	-11/28&	-5/14&	-9/28&	-2/7&	-1/4&	-3/14&	-5/28&	-1/7&	-3/28&	-1/14&	-1/28\\%
&
0/1&	1/28&	1/14&	3/28&	1/7&	5/28&	3/14&	1/4&	2/7&	9/28&	5/14&	11/28&	3/7&	13/28\\
\bottomrule
\end{tabular}
\caption[font=tiny,format=hang,justification=left]{$|r|=1$, $s \equiv 0, p_1 \equiv 0 \bmod r$.} \label{ta:hc_1_m0_m0}
\end{sidewaystable}
}

\begin{sidewaystable}
\newcolumntype{C}{>{\tiny $}r<{$}}
\begin{tabular}{*{15}{C}} \toprule
\multicolumn{15}{c}{Computed Kreck-Stolz invariants of \eks{s}}\\

\multicolumn{15}{c}{with $|r|=3$, $s \equiv -1, p_1 \equiv 0 \bmod r$.} \\
\multicolumn{1}{c}{$s_2$} & \multicolumn{14}{c}{$s_1$} \\\midrule
\multirow{2}{8mm}{$-5/18$} &-55/112&	-51/112&	-47/112&	-43/112&	-39/112&	-5/16&	-31/112&	-27/112&	-23/112&	-19/112&	-15/112&	-11/112&	-1/16&	-3/112\\%
&
1/112&	5/112&	9/112&	13/112&	17/112&	3/16&	25/112&	29/112&	33/112&	37/112&	41/112&	45/112&	7/16&	53/112\\

\multirow{2}{8mm}{$-1/36$} &-55/112&	-51/112&	-47/112&	-43/112&	-39/112&	-5/16&	-31/112&	-27/112&	-23/112&	-19/112&	-15/112&	-11/112&	-1/16&	-3/112\\%
&
1/112&	5/112&	9/112&	13/112&	17/112&	3/16&	25/112&	29/112&	33/112&	37/112&	41/112&	45/112&	7/16&	53/112\\

\multirow{2}{8mm}{$2/9$} &-55/112&	-51/112&	-47/112&	-43/112&	-39/112&	-5/16&	-31/112&	-27/112&	-23/112&	-19/112&	-15/112&	-11/112&	-1/16&	-3/112\\
\bottomrule
\end{tabular}
\caption{$|r|=3$, $s \equiv \pm1, p_1 \equiv 0 \bmod r$.} \label{ta:hc_3_p1_m0}
\end{sidewaystable}

}
\end{mytable}

\begin{mytable}
{\tiny
\newcolumntype{R}{>{$}r<{$}}
\begin{longtable}[c]{*{2}{*{5}{R}}}
\caption{A summary of the invariants of \eks{s}, where $k \in
  \Z_{12}$, $l \in \Z_{28}$ and $s_1,s_2$ take values in
  $\left[-\frac{1}{2},\frac{1}{2}\right) \bmod 1$.} \label{ta:eks-invariants}
\endfirsthead
\toprule
\multicolumn{10}{c}{\tiny\slshape Table \ref{ta:eks-invariants}, continued from previous page} \\
\toprule
\absr & s & p_1 & s_2 & s_1 & \absr & s & p_1 & s_2 & s_1 \\
\midrule
\endhead
\bottomrule
\multicolumn{10}{c}{\tiny\slshape continued next page} \\
\endfoot
\bottomrule
\endlastfoot
\absr & s & p_1 & s_2 & s_1 & \absr & s & p_1 & s_2 & s_1 \\
\midrule
1	& 0	& 0	& 0/1+k/12	& -1/2+l/28&
&&&& \\
3	& -1	& 0	& 2/9+k/4	& -55/112+l/28&
3	& 1	& 0	& 1/36+k/4	& -53/112+l/28\\
5	& -2	& 2	& 0/1+k/12	& -69/140+l/28&
5	& -1	& -2	& 1/60+k/12	& -131/280+l/28\\
5	& 1	& -2	& 1/15+k/12	& -139/280+l/28&
5	& 2	& 2	& 0/1+k/12	& -33/70+l/28\\
7	& -3	& -3	& 1/84+k/12	& -389/784+l/28&
7	& -3	& 1	& 1/28+k/12	& -373/784+l/28\\
7	& -3	& 2	& 0/1+k/12	& -365/784+l/28&
7	& -2	& 0	& 1/42+k/12	& -53/112+l/28\\
7	& -1	& 0	& 1/21+k/12	& -53/112+l/28&
7	& 1	& 0	& 1/28+k/12	& -55/112+l/28\\
\end{longtable}
}
\end{mytable}

\subsection{Appendix C} \label{ssec:app-c}
Tables \ref{ta:eks-r-1}--\ref{ta:eks-r-5} list homeomorphism
classes of \eks{s} with the rank of the fourth integral cohomology
group equal to $|r|=1,3,5$ respectively. Each smooth structure in each
such homeomorphism class is represented by an \eks{}; these tables
list the `smallest' representatives. 

\begin{mytable}
{\tiny
\newcolumntype{R}{>{$}r<{$}}
\begin{longtable}[c]{RR}
\caption{Smallest \eks{s} with $|r|=1$, $s_2 \geq 0$.} \label{ta:eks-r-1}
\endfirsthead
\multicolumn{2}{c}{\tiny\slshape Table \ref{ta:eks-r-1}, continued from previous page} \\
\endhead
\multicolumn{2}{c}{\tiny\slshape continued next page} \\
\endfoot
\endlastfoot
\begin{tabular}{*{6}{R}}
\toprule
\multicolumn{6}{c}{\minhcheadera{}} \\
\multicolumn{6}{c}{\minhcheaderb{}} \\
\multicolumn{6}{c}{\minhcheaderc{1}{0}{0}{0/1}} \\
{\rm sum} & k_0 & k_1 & l_0 &l_1 & s_1 \\
\midrule
0	& -17	& 5	& -14	& -2	& -1/2\\
0	& -24	& 6	& -13	& -12	& -13/28\\
0	& -47	& 11	& -38	& -8	& -3/7\\
0	& -7	& 3	& -6	& 0	& -11/28\\
0	& -17	& 5	& -16	& 2	& -5/14\\
0	& -36	& 18	& -29	& -4	& -9/28\\
0	& -27	& -1	& -24	& -6	& -2/7\\
0	& 7	& 5	& 6	& 6	& -1/4\\
0	& -41	& -3	& -38	& -8	& -3/14\\
0	& -12	& -1	& -10	& -4	& -5/28\\
0	& -17	& -4	& -14	& -8	& -1/7\\
0	& -11	& -9	& -10	& -10	& -3/28\\
0	& -24	& 6	& -23	& 3	& -1/14\\
0	& 3	& 1	& 2	& 2	& -1/28\\
0	& -4	& 1	& -2	& -2	& 0/1\\
0	& -3	& -1	& -2	& -2	& 1/28\\
0	& -20	& -3	& -18	& -6	& 1/14\\
0	& 11	& 9	& 10	& 10	& 3/28\\
0	& -22	& 8	& -21	& 4	& 1/7\\
0	& -14	& 4	& -13	& 1	& 5/28\\
0	& -46	& 8	& -44	& 3	& 3/14\\
0	& -7	& -5	& -6	& -6	& 1/4\\
0	& -30	& 6	& -28	& 1	& 2/7\\
0	& -33	& 4	& -18	& -18	& 9/28\\
0	& -14	& -2	& -12	& -5	& 5/14\\
0	& -6	& 0	& -4	& -3	& 11/28\\
0	& -46	& 8	& -36	& -11	& 3/7\\
0	& -25	& 12	& -18	& -6	& 13/28\\

\bottomrule
\end{tabular}
&
\begin{tabular}{*{6}{R}}
\toprule
\multicolumn{6}{c}{\minhcheadera{}} \\
\multicolumn{6}{c}{\minhcheaderb{}} \\
\multicolumn{6}{c}{\minhcheaderc{1}{0}{0}{1/12}} \\
{\rm sum} & k_0 & k_1 & l_0 &l_1 & s_1 \\
\midrule
1	& -26	& 6	& -23	& -1	& -1/2\\
-2	& 10	& 8	& 9	& 9	& -13/28\\
1	& -14	& 5	& -7	& -7	& -3/7\\
1	& -21	& -3	& -19	& -6	& -11/28\\
1	& -11	& 4	& -9	& -1	& -5/14\\
1	& -29	& 9	& -22	& -6	& -9/28\\
1	& -38	& 2	& -33	& -7	& -2/7\\
1	& -30	& 10	& -19	& -11	& -1/4\\
1	& -13	& -1	& -11	& -4	& -3/14\\
1	& -6	& 2	& -3	& -3	& -5/28\\
1	& -14	& -3	& -11	& -7	& -1/7\\
1	& -13	& 5	& -10	& -2	& -3/28\\
1	& -34	& -2	& -25	& -15	& -1/14\\
1	& -17	& -1	& -12	& -8	& -1/28\\
1	& -43	& 19	& -42	& 13	& 0/1\\
1	& -17	& -5	& -12	& -11	& 1/28\\
1	& -11	& -2	& -7	& -7	& 1/14\\
1	& -10	& 5	& -9	& 1	& 3/28\\
1	& -27	& 8	& -21	& -5	& 1/7\\
1	& -42	& 14	& -27	& -15	& 5/28\\
1	& -21	& 11	& -20	& 5	& 3/14\\
1	& -3	& -1	& -2	& -2	& 1/4\\
1	& -17	& 5	& -11	& -6	& 2/7\\
1	& -43	& 5	& -34	& -11	& 9/28\\
1	& -20	& 0	& -17	& -5	& 5/14\\
-2	& 2	& 0	& 1	& 1	& 11/28\\
1	& -5	& 3	& -4	& 0	& 3/7\\
1	& -13	& 7	& -11	& 0	& 13/28\\

\bottomrule
\end{tabular}
\\
\end{longtable}

\begin{longtable}[c]{RR}
\caption{Smallest \eks{s} with $|r|=3$, $s_2 \geq 0$.} \label{ta:eks-r-3}
\endfirsthead
\multicolumn{2}{c}{\tiny\slshape Table \ref{ta:eks-r-3}, continued from previous page} \\
\endhead
\multicolumn{2}{c}{\tiny\slshape continued next page} \\
\endfoot
\endlastfoot
\begin{tabular}{*{6}{R}}
\toprule
\multicolumn{6}{c}{\minhcheadera{}} \\
\multicolumn{6}{c}{\minhcheaderb{}} \\
\multicolumn{6}{c}{\minhcheaderc{3}{-1}{0}{2/9}} \\
{\rm sum} & k_0 & k_1 & l_0 &l_1 & s_1 \\
\midrule
0	& -165	& 60	& -152	& 16	& -55/112\\
0	& -129	& 45	& -128	& 40	& -51/112\\
0	& -63	& -36	& -56	& -44	& -47/112\\
0	& -51	& 15	& -32	& -20	& -43/112\\
0	& -132	& 27	& -128	& 16	& -39/112\\
0	& -164	& -20	& -156	& -33	& -5/16\\
0	& -36	& -15	& -32	& -20	& -31/112\\
0	& -68	& 16	& -60	& -3	& -27/112\\
0	& -39	& 12	& -20	& -20	& -23/112\\
0	& -156	& 45	& -92	& -68	& -19/112\\
0	& -104	& -20	& -81	& -51	& -15/112\\
0	& -164	& 52	& -108	& -57	& -11/112\\
0	& -132	& 51	& -104	& -20	& -1/16\\
0	& -81	& 36	& -80	& 28	& -3/112\\
0	& -12	& -3	& -8	& -8	& 1/112\\
0	& -60	& -27	& -56	& -32	& 5/112\\
0	& -177	& 36	& -164	& 4	& 9/112\\
0	& -60	& 15	& -56	& 4	& 13/112\\
0	& -57	& 21	& -56	& 16	& 17/112\\
0	& -92	& 40	& -75	& -9	& 3/16\\
0	& -84	& 9	& -68	& -20	& 25/112\\
0	& -56	& 16	& -36	& -21	& 29/112\\
0	& -75	& -9	& -68	& -20	& 33/112\\
0	& -33	& -3	& -20	& -20	& 37/112\\
0	& -21	& 9	& -20	& 4	& 41/112\\
0	& -15	& 3	& -8	& -8	& 45/112\\
0	& -84	& 33	& -80	& 16	& 7/16\\
0	& -129	& 36	& -104	& -20	& 53/112\\

\bottomrule
\end{tabular}
&
\begin{tabular}{*{6}{R}}
\toprule
\multicolumn{6}{c}{\minhcheadera{}} \\
\multicolumn{6}{c}{\minhcheaderb{}} \\
\multicolumn{6}{c}{\minhcheaderc{3}{-1}{0}{17/36}} \\
{\rm sum} & k_0 & k_1 & l_0 &l_1 & s_1 \\
\midrule
0	& -2	& 1	& 0	& 0	& -55/112\\
0	& -66	& 12	& -41	& -29	& -51/112\\
0	& -89	& 19	& -84	& 6	& -47/112\\
0	& -29	& 7	& -18	& -12	& -43/112\\
0	& -35	& -11	& -24	& -24	& -39/112\\
0	& -50	& 7	& -30	& -24	& -5/16\\
0	& -6	& 0	& -5	& -2	& -31/112\\
0	& -42	& 12	& -26	& -17	& -27/112\\
0	& -30	& 12	& -29	& 7	& -23/112\\
0	& -54	& 24	& -50	& 7	& -19/112\\
0	& -30	& 0	& -29	& -2	& -15/112\\
0	& -174	& 60	& -170	& 43	& -11/112\\
0	& -24	& 0	& -23	& -2	& -1/16\\
0	& -66	& 6	& -53	& -17	& -3/112\\
0	& -42	& 0	& -41	& -2	& 1/112\\
0	& -54	& 18	& -29	& -26	& 5/112\\
0	& -41	& -2	& -30	& -18	& 9/112\\
0	& -36	& -18	& -29	& -26	& 13/112\\
0	& -90	& 42	& -77	& -2	& 17/112\\
0	& -114	& -12	& -74	& -65	& 3/16\\
0	& -36	& 0	& -35	& -2	& 25/112\\
0	& -30	& -12	& -26	& -17	& 29/112\\
0	& -72	& 6	& -50	& -29	& 33/112\\
0	& -48	& 18	& -41	& -2	& 37/112\\
0	& -12	& 0	& -11	& -2	& 41/112\\
0	& -65	& 19	& -42	& -24	& 45/112\\
0	& -18	& 0	& -17	& -2	& 7/16\\
0	& -114	& 30	& -65	& -53	& 53/112\\

\bottomrule
\end{tabular}
\\
\end{longtable}

\begin{longtable}[c]{RR}
\caption{Smallest \eks{s} with $|r|=5$, $s_2 \geq 0$.} \label{ta:eks-r-5}
\endfirsthead
\multicolumn{2}{c}{\tiny\slshape Table \ref{ta:eks-r-5}, continued from previous page} \\
\endhead
\multicolumn{2}{c}{\tiny\slshape continued next page} \\
\endfoot
\endlastfoot
\begin{tabular}{*{6}{R}}
\toprule
\multicolumn{6}{c}{\minhcheadera{}} \\
\multicolumn{6}{c}{\minhcheaderb{}} \\
\multicolumn{6}{c}{\minhcheaderc{5}{-2}{2}{0/1}} \\
{\rm sum} & k_0 & k_1 & l_0 &l_1 & s_1 \\
\midrule
0	& -168	& 52	& -161	& 28	& -69/140\\
0	& -128	& 12	& -92	& -47	& -16/35\\
0	& -81	& -27	& -64	& -48	& -59/140\\
0	& -228	& 112	& -177	& -36	& -27/70\\
0	& -76	& -7	& -48	& -44	& -7/20\\
0	& -87	& 3	& -68	& -28	& -11/35\\
0	& -124	& -48	& -121	& -52	& -39/140\\
0	& -187	& -1	& -144	& -68	& -17/70\\
0	& -88	& 16	& -57	& -36	& -29/140\\
0	& -147	& 39	& -84	& -68	& -6/35\\
0	& -116	& 13	& -108	& -4	& -19/140\\
0	& -47	& 4	& -28	& -24	& -1/10\\
0	& -168	& 16	& -157	& -7	& -9/140\\
0	& -161	& 28	& -88	& -84	& -1/35\\
0	& -228	& 56	& -157	& -76	& 1/140\\
0	& -164	& 72	& -156	& 33	& 3/70\\
0	& -12	& 3	& -8	& -4	& 11/140\\
0	& -172	& 83	& -168	& 52	& 4/35\\
0	& -52	& -1	& -48	& -8	& 3/20\\
0	& -64	& 12	& -37	& -31	& 13/70\\
0	& -108	& 32	& -107	& 28	& 31/140\\
0	& -188	& 56	& -101	& -92	& 9/35\\
0	& -124	& -8	& -77	& -71	& 41/140\\
0	& -67	& 23	& -64	& 12	& 23/70\\
0	& -127	& 59	& -108	& -4	& 51/140\\
0	& -212	& 103	& -204	& 52	& 2/5\\
0	& -144	& 32	& -137	& 13	& 61/140\\
0	& -24	& -8	& -21	& -12	& 33/70\\

\bottomrule
\end{tabular}
&
\begin{tabular}{*{6}{R}}
\toprule
\multicolumn{6}{c}{\minhcheadera{}} \\
\multicolumn{6}{c}{\minhcheaderb{}} \\
\multicolumn{6}{c}{\minhcheaderc{5}{-2}{2}{1/12}} \\
{\rm sum} & k_0 & k_1 & l_0 &l_1 & s_1 \\
\midrule
1	& -69	& 35	& -42	& -26	& -69/140\\
1	& -96	& 44	& -95	& 35	& -16/35\\
1	& -99	& 25	& -91	& 4	& -59/140\\
1	& -35	& 15	& -32	& 4	& -27/70\\
1	& -86	& 18	& -85	& 15	& -7/20\\
1	& -55	& 25	& -51	& 8	& -11/35\\
1	& -59	& -5	& -52	& -16	& -39/140\\
1	& -95	& -5	& -91	& -12	& -17/70\\
1	& -27	& -6	& -25	& -9	& -29/140\\
1	& -75	& 5	& -51	& -32	& -6/35\\
1	& -45	& 21	& -32	& -11	& -19/140\\
1	& -106	& 53	& -89	& -5	& -1/10\\
1	& -91	& 24	& -79	& -5	& -9/140\\
1	& -91	& 18	& -65	& -29	& -1/35\\
1	& -91	& 44	& -55	& -35	& 1/140\\
1	& -36	& 8	& -35	& 5	& 3/70\\
1	& -65	& 11	& -46	& -22	& 11/140\\
1	& -75	& 31	& -67	& 4	& 4/35\\
1	& -25	& 1	& -16	& -12	& 3/20\\
1	& -109	& 55	& -91	& -6	& 13/70\\
1	& -149	& 75	& -86	& -62	& 31/140\\
1	& -62	& -6	& -49	& -25	& 9/35\\
1	& -11	& -2	& -9	& -5	& 41/140\\
1	& -96	& 13	& -75	& -25	& 23/70\\
1	& -32	& -11	& -29	& -15	& 51/140\\
1	& -16	& 4	& -15	& 1	& 2/5\\
1	& -55	& 15	& -51	& 4	& 61/140\\
1	& -105	& 21	& -76	& -32	& 33/70\\

\bottomrule
\end{tabular}
\\
\end{longtable}

\begin{extratable}
\begin{longtable}[c]{RR}
\caption{Smallest \eks{s} with $|r|=7$, $s_2 \geq 0$.} \label{ta:eks-r-7}
\endfirsthead
\multicolumn{2}{c}{\tiny\slshape Table \ref{ta:eks-r-5}, continued from previous page} \\
\endhead
\multicolumn{2}{c}{\tiny\slshape continued next page} \\
\endfoot
\endlastfoot
\begin{tabular}{|*{6}{R}|}
\hline
\multicolumn{6}{|c|}{\minhcheadera{}} \\
\multicolumn{6}{|c|}{\minhcheaderb{}} \\
\multicolumn{6}{|c|}{\minhcheaderc{7}{-3}{-3}{1/84}} \\
\hline
{\rm sum} & k_0 & k_1 & l_0 &l_1 & s_1 \\
\hline\hline
1	& -92	& -4	& -63	& -45	& -389/784\\
1	& -77	& 39	& -71	& 10	& -361/784\\
1	& -171	& 77	& -144	& -8	& -333/784\\
1	& -106	& -39	& -77	& -73	& -305/784\\
1	& -204	& 97	& -175	& -3	& -277/784\\
1	& -288	& 48	& -283	& 35	& -249/784\\
1	& -92	& 20	& -87	& 7	& -221/784\\
1	& -172	& 76	& -171	& 67	& -193/784\\
1	& -73	& 21	& -39	& -36	& -165/784\\
1	& -36	& 17	& -35	& 11	& -137/784\\
1	& -95	& -8	& -87	& -21	& -109/784\\
1	& -185	& 25	& -183	& 20	& -81/784\\
1	& -144	& 41	& -119	& -17	& -53/784\\
1	& -287	& 67	& -172	& -127	& -25/784\\
1	& -71	& -18	& -49	& -45	& 3/784\\
1	& -101	& 21	& -88	& -8	& 31/784\\
1	& -64	& 24	& -59	& 7	& 59/784\\
1	& -127	& 20	& -77	& -59	& 87/784\\
1	& -36	& -15	& -31	& -21	& 115/784\\
1	& -64	& 17	& -35	& -31	& 143/784\\
1	& -207	& 104	& -203	& 67	& 171/784\\
1	& -143	& 35	& -127	& -4	& 199/784\\
1	& -120	& 24	& -115	& 11	& 227/784\\
1	& -35	& -3	& -32	& -8	& 255/784\\
1	& -95	& 48	& -87	& 11	& 283/784\\
1	& -183	& 48	& -115	& -73	& 311/784\\
1	& -204	& 76	& -199	& 53	& 339/784\\
1	& -87	& 39	& -71	& -8	& 367/784\\

\hline
\end{tabular}
&
\begin{tabular}{|*{6}{R}|}
\hline
\multicolumn{6}{|c|}{\minhcheadera{}} \\
\multicolumn{6}{|c|}{\minhcheaderb{}} \\
\multicolumn{6}{|c|}{\minhcheaderc{7}{-3}{-3}{2/21}} \\
\hline
{\rm sum} & k_0 & k_1 & l_0 &l_1 & s_1 \\
\hline\hline
-1	& 60	& -6	& 35	& 31	& -389/784\\
-1	& 81	& 50	& 73	& 59	& -361/784\\
-1	& 190	& -31	& 175	& 3	& -333/784\\
-1	& 8	& 1	& 7	& 3	& -305/784\\
-1	& 119	& -53	& 113	& -24	& -277/784\\
-1	& 137	& -34	& 105	& 31	& -249/784\\
-1	& 78	& -24	& 59	& 17	& -221/784\\
-1	& 4	& 1	& 3	& 3	& -193/784\\
-1	& 87	& -25	& 64	& 22	& -165/784\\
-1	& 143	& -11	& 81	& 78	& -137/784\\
-1	& 59	& -11	& 57	& -6	& -109/784\\
-1	& 18	& 1	& 17	& 3	& -81/784\\
-1	& 130	& 1	& 87	& 63	& -53/784\\
-1	& 81	& -20	& 63	& 17	& -25/784\\
-1	& 283	& -105	& 281	& -94	& 3/784\\
-1	& 32	& 1	& 31	& 3	& 31/784\\
-1	& 143	& -67	& 134	& -24	& 59/784\\
-1	& 77	& -25	& 57	& 18	& 87/784\\
-1	& 105	& 59	& 88	& 78	& 115/784\\
-1	& 46	& 1	& 45	& 3	& 143/784\\
-1	& 172	& -62	& 87	& 87	& 171/784\\
-1	& 73	& 35	& 60	& 50	& 199/784\\
-1	& 106	& 57	& 101	& 63	& 227/784\\
-1	& 36	& 1	& 35	& 3	& 255/784\\
-1	& 49	& 3	& 46	& 8	& 283/784\\
-1	& 190	& -24	& 119	& 87	& 311/784\\
-1	& 92	& 1	& 73	& 31	& 339/784\\
-1	& 22	& 1	& 21	& 3	& 367/784\\

\hline
\end{tabular}
\\
\end{longtable}
\end{extratable}
}
\end{mytable}

\subsection{Appendix D} \label{ssec:app-d}
The graphs in figure \ref{gr:tt} depict the outcome of some time
trials that compared alternative methods of computing the invariants
of the \eks{s}. The \cpp{} code is significantly faster. The sole
exception to this statement is when $\max |k| \leq 4$. The reason for
this apparent anomaly is that the \cpp{} code search and found the
\eks{s} in addition to computing their invariants; the \maple{} code
merely computed the invariants.

{
\def\timesec{{\begin{rotate}{90}Time (secs)\end{rotate}}}
\def\timehr{{\begin{rotate}{90}Time (hrs)\end{rotate}}}
\def\Diamond{{{\rm x}}}
\def\Box{{*}}
\def\cppx{{\rm x}}
\def\maplex{{*}}
\begin{center}
\begin{figure}[!htb]
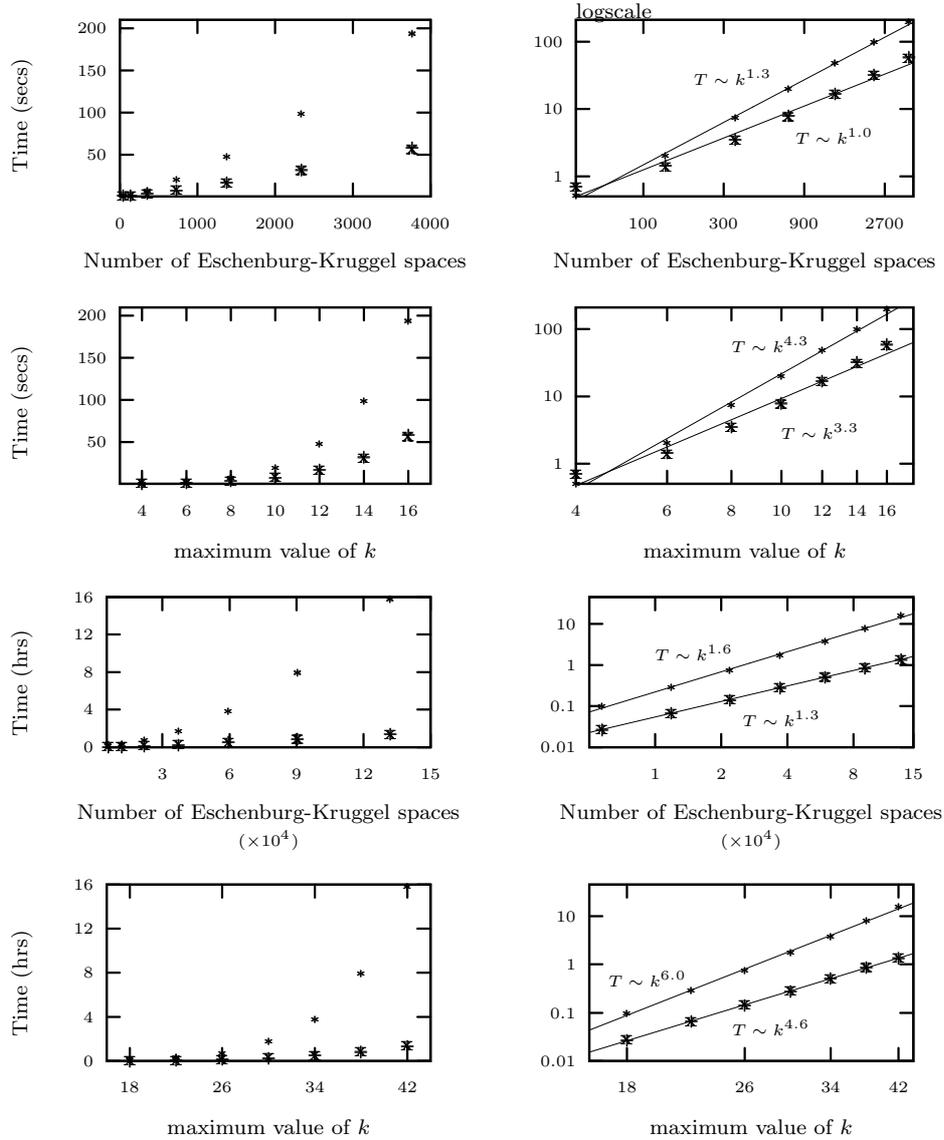

\input{time_trial_2_16.eepic}
\input{time_trial_18_42.eepic}
\caption{Even rows: Computation time versus the number, $N=N(k)$, of
  Eschenburg-Kruggel spaces in the cube $[-k,k]^6$ with $\sum_i
  k_i=\sum_i l_i \in [0,2]$; Odd rows: Computation time versus the
  number $k$. The right column shows the same data in
  log-scale. $*$=\maple{}~\cite{maple} times;
  $+,\Diamond$=\cpp-times.} \label{gr:tt}
\end{figure}
\end{center}
}

The time trials were conducted on a quad-core Intel Xeon 5148 with an
over-clocked 2.33GHz cpu, a 4MB cache, 3.2GB ram memory, and 8.4GB swap
memory.



\end{document}